\documentclass[journal]{IEEEtran}
\usepackage{blindtext}
\usepackage{graphicx}
\usepackage{amssymb}
\usepackage{amsmath}
\usepackage{amsthm}
\usepackage{setspace}
\usepackage{amsthm}
\usepackage{mathrsfs}
\usepackage{upgreek}
\usepackage{algorithm2e}
\usepackage{multirow}
\usepackage{multicol}
\usepackage{float}

\usepackage[hidelinks]{hyperref}
\usepackage{subfig}
\usepackage{bm}
\usepackage{color}

\usepackage[letterpaper, margin=1in]{geometry}

\newcommand*\tabrule[1][]{%
   \if\relax\detokenize{#1}\relax
     \rule{\linewidth}{\arrayrulewidth}%
   \else
     \rule{#1}{\arrayrulewidth}%
   \fi
}

\newtheorem{theorem}{Theorem}[section]

\newtheorem{assumption}[theorem]{Assumption}

\usepackage{cite}

\ifCLASSINFOpdf
\else
\fi
\begin{document}

\singlespacing

\title{Optimization Models for Estimating Transit Network Origin-Destination Flows with AVL/APC Data}

\author{Xinyu Liu, 
        Pascal~Van~Hentenryck,
        and Xilei Zhao
\thanks{X. Liu and P. Van Hentenryck are with the H. Milton Stewart School of Industrial and Systems Engineering, Georgia Institute of Technology, Atlanta, GA, 30332 USA.}
\thanks{X. Zhao is with the Department of Civil and Coastal Engineering, University of Florida, Gainesville, FL, 32611 USA (Corresponding author: Xilei Zhao, e-mail: xilei.zhao@essie.ufl.edu)}}

\maketitle

\begin{abstract}

With the increasing adoption of Automatic Vehicle Location (AVL) and Automatic Passenger Count (APC) technologies by transit agencies, a massive amount of time-stamped and location-based passenger boarding and alighting count data can be collected on a continuous basis. The availability of such large-scale transit data offers new opportunities to produce estimates for Origin-Destination (O-D) flows, helping inform transportation planning and transit management. However, the state-of-the-art methodologies for AVL/APC data analysis mostly tackle the O-D flow estimation problem within routes and barely infer the transfer activities across the entire transit network. This paper proposes three optimization models to identify transfers and approximate network-level O-D flows by minimizing the deviations between estimated and observed proportions or counts of transferring passengers: A Quadratic Integer Program (QIP), a feasible rounding procedure for the Quadratic Convex Programming (QCP) relaxation of the QIP, and an Integer Program (IP). The inputs of the models are readily available by applying the various route-level flow estimation algorithms to the automatically collected AVL/APC data and the output of the models is a network O-D estimation at varying geographical resolutions. The optimization models were evaluated on a case study for Ann Arbor-Ypsilanti area in Michigan. The IP model outperforms the QCP approach in terms of accuracy and remains tractable from an efficiency standpoint, contrary to the QIP. Its estimated O-D matrix achieves an R-Squared metric of 95.57\% at the Traffic Analysis Zone level and 92.39\% at the stop level, compared to the ground-truth estimates inferred from the state-of-practice trip-chaining methods.

\end{abstract}


%
\IEEEpeerreviewmaketitle

\section{Introduction} 

The transit origin-destination (O-D) matrix is a major input for public transit agencies to conduct scheduling and operations planning, long-term planning and design, performance analysis, and market evaluations. However, the traditional data sources for transit O-D matrices are transit on-board surveys, a time-consuming and labor-intensive collection process that is prone to sampling errors. Fortunately, recent years have witnessed a growing interest in building transit O-D matrices from data sources that are automatically collected through intelligent transportation systems, such as Automatic Vehicle Location (AVL), Automatic Passenger Count (APC), and Automatic Fare Collection (AFC) systems \cite{Iliopoulou2019ITStransportplanning}. 

In particular, automatic passenger counter is an electronic device usually installed on transit vehicles which records boarding and alighting data. Enabled by technologies such as infrared or light beams, digital cameras, thermal imaging and ultrasonic detection, the collected data is of high accuracy and can be easily validated \cite{Wiki2019APC,ITS2019APC}. APC technology is commonly deployed with AVL technology which provides access to real-time transit vehicle dispatching and tracking data through information technology and Global Positioning Systems \cite{Wiki2019AVL}. Most transit agencies have installed APC systems on at least 10\% to 15\% of their bus fleet \cite{Furth2005APCmainstream} and AVL is expected to be present in most fixed-route systems \cite{Wiki2019AVL}, as well as bus rapid transit systems \cite{Parker2008AVLUpdates}. The abundant AVL/APC data can jointly link passenger data to vehicle location \cite{Wiki2019APC} and thus offers a rich source of data in both spatial and temporal dimensions. 

Ever since the proliferation of AVL/APC technologies, researchers have worked on the conceptualization and development of methodologies to exploit the large-scale transit data they collect for use in transit performance analysis and travel demand modeling. Traditionally, such data contributed greatly to transit performance analysis and service management: They mainly address the problems of determining vehicle loads or run times \cite{Tetreault2010RunTime}, diagnosing or improving transit system performance \cite{Mandelzys2010AVLAPC, Furth2003AVLAPC}, and analyzing transit ridership \cite{Furth2003AVLAPC, Golani2007AVLAPC}. 
Recently, many efforts have been directed to O-D flow estimation based on AVL/APC data. Iterative Proportion Filtering (IPF) procedure is one of the most widely accepted methodologies: It aims at estimating population O-D flow on each transit route based on sampled stop-level boarding and alighting counts. Deming and Stephan \cite{Deming1940IPF} first proposed a procedure to adjust sampled frequency tables with known marginal totals obtained from different sources. The method is a good fit for route-level O-D estimation where a base matrix sampled from on-board surveys is adjusted by boarding and alighting counts at each stop. Ben-Akiva \cite{BenAkiva1985firstIPF} then showed IPF to be cost-effective for route-level O-D table estimation when combined with on-board survey data. However, IPF has limitations for its dependence on a base matrix constructed from on-board surveys. The vast amount of AVL/APC data made available recently reduces this dependence in many ways, including by enabling more cost-effective choices of base matrices and by inspiring new methods that require little survey data. Using the APC data available by the campus transit service in Ohio State University, McCord et al. \cite{McCord2010routeIPF} demonstrated that the IPF procedure, even with a non-informative, null base matrix, can achieve comparable O-D estimates as on-board surveys. Empirical results then showed that assuming no \textit{a priori} estimate of O-D flows, an arbitrary base matrix that is adjusted iteratively with APC data can achieve a higher accuracy than the null base matrix \cite{Ji2014IterativeIPF}. Moreover, Ji et al. \cite{Ji2015Survey&APC} developed a heuristic expectation-maximization method using APC and on-board survey data, which was shown to outperform the IPF procedure when little survey data is present. In \cite{Ji2015StatsRouteCounts}, a Markov chain Monte Carlo simulation approach with sampling is proposed to infer route O-D with large amounts of APC data only and validated through a numerical test. These advances in exploiting AVL/APC data allow for accurate extrapolation of O-D flows along the specific route even in absence of costly, time-consuming and error-prone survey data.

However, it is important to realize that the route-level O-D matrix, though helpful to transit planning in terms of designing route patterns and service frequencies \cite{McCord2010routeIPF}, does not reflect the true O-D flows of transit passengers. The route-level O-D matrix only represents the flow distribution along a single route, while the true passenger trajectories might include additional travels to and from outside this route \cite{Ji2015StatsRouteCounts}. Thus there is an identified need to infer transit network O-D flow matrices based on AVL/APC data. However, the state-of-the-art analysis techniques targeting transit network level is still lacking. As APC data only records the time-stamped number of stop-level boardings and alightings for each route, they cannot differentiate initial and transfer boardings \cite{chu2004ridership}. Also, most route-level O-D estimation methods cannot generalize to transit network if only AVL/APC data is available \cite{McCord2010routeIPF,Ji2015Survey&APC,Ji2015StatsRouteCounts}. Thus a transfer identification algorithm is required to resolve this issue and to produce an overall O-D matrix on the transit network. 

The goal of this paper is thus to propose novel optimization models to identify transfer activities from AVL/APC data and observed proportion of transferring passengers at different transit centers and other stops. 
As a result, the paper addresses a limitation of methodologies for APC data analysis, by extending the state-of-the-art route-level O-D matrix estimation and producing a network-level O-D matrix. Observe that the O-D matrix is not intended to predict or analyze individual activities, but rather to assess travel behavior and demand at an aggregate level. This is perfectly reasoned given that the aggregate nature of APC data already makes it challenging to accurately recover the path choices of each rider \cite{Ji2015Survey&APC, Ji2015StatsRouteCounts}. Moreover, many applications in transit planning do not need individual-level activities, but need aggregate-level travel demand information, e.g., at the Traffic Analysis Zone (TAZ) level or other self-defined geographical boundaries. The optimization models proposed in this paper generate aggregate O-D matrices at different geographical resolutions and hence can inform future transit planning and investment decisions. 

The rest of the paper is organized as follows. The next section defines the problem with a specification of available input data and desired output data. The following section proposes three optimization models to solve the transfer identification problem. These models are then evaluated on a case study. The paper concludes by summarizing the results and discussing its practical applications in transportation planning.

\section{Problem Definition}

The time-stamped route O-D flow information is made accessible for transportation planners and researchers through wide adoption of APC technologies for raw data collection and well-developed methodologies for estimating population route O-D flows \cite{McCord2010routeIPF, Ji2014IterativeIPF, Ji2015Survey&APC, Ji2015StatsRouteCounts}. The obtained boarding-alighting pairs or route-level O-D flows do not equate to O-D transit trips \cite{McCord2010routeIPF}. Indeed, due to the transfer activities, many actual trips may contain several boarding-alighting pairs where each is called a trip segment or trip leg. A trip can be generally represented as one segment or a sequence of segment(s), where the former is referred to as a singleton trip in this study. Each trip segment in the latter is described by a corresponding ordinal number; for example, a trip with two segments has a first leg and a second leg. Accordingly, there is a need for transfer identification when integrating O-D flows on different routes to generate the transit network O-D matrix. This Transfer Identification Problem (TIP) is the focus of this paper and the rest of this section specifies its input and output.

\subsection{Input Specification}

Since transfers are based on individual activities, identifying them first requires disaggregation of time-stamped route O-D flows into individual records each specifying the route, the boarding and alighting times and stops. Equivalently, it requires the time-stamped passenger counts describing the demand at a stop as trip origin or destination, and the route O-D matrices specifying the distribution of alighting stops for all boarding passengers at a stop. The model also needs access to the bus schedules for each route and the stop locations on all routes. In addition, the model assumes the availability of observed or estimated proportions of transfers\footnote{The terms transfer probabilities and transfer rates are referred to as the proportions of transfers.}. Such data can be obtained directly from relatively long-term observations at stop terminals or estimated by transportation professionals based on years of experience. This data collection process outperforms on-board surveys in terms of ease to implement and reliability: Enumerating the passengers leaving the terminals requires less efforts compared to conducting detailed surveys for individual passengers, and it is more reliable than on-board surveys by reflecting the whole population rather than a proportion of passengers on selected routes. Moreover, the observed transfer probabilities can be obtained by averaging over a large amount of historical observations, thus less prone to inaccurate or misrecorded data entries on individual passengers.

\subsection{Transfer Assumptions}
The transfer activities modeled in this work are characterized by three behavioral assumptions as stated, justified and discussed below.

\begin{assumption}
    Transfer activities that happen in transit centers should be differentiated from those occurring at non-transit centers.
\end{assumption}

It follows from the observation that more transfer trips are expected at transit centers, which are designed to be served by multiple bus or rail routes synchronized for facilitating transfers. Consequently, the observed transfer probabilities at transit centers are expected to be more significant than at other stops. This model thus evaluates them separately. 

\begin{assumption}
A transfer between two trip segments is only feasible when the following three conditions are satisfied. First, the two trip segments linked by a transfer must not belong to the same route. Second, passengers only transfer within some thresholds for walking distance and transfer time. Third, transfers are directional such that the prior trip segment must have ended earlier than the boarding of the subsequent one. 
\end{assumption}

These assumptions on travel behaviors are widely applied and tested for generating trip chains \cite{Alsger2015SmartCard, trepanier2007individual, barry2002origin, munizaga2012estimation}. To verify the feasibility of a transfer between any two trips in terms of walking distance and transfer time intervals, the stop locations and bus schedules are required, which are often made publicly available from local regulators. In particular, this information can be extracted from the General Transit Feed Specification (GTFS) data, which is a common format for public transportation schedules and the associated geographic information \cite{GoogleTransitAPIs2019}. Notably, the route and stop ID information can jointly map each bus stop to a geographical location; the alighting times can be calculated from boarding times and bus schedules.

\begin{assumption}
    Passengers transfer at most once.
\end{assumption}

In general, most trips in a bus-based transit network involve no transfer or a single transfer. This study adopts this \textit{one-transfer assumption}. The model however can easily handle the more general case (i.e., two or more transfers) if needed.

\subsection{Output Specification}

The desired output for this study is a transit network level O-D flow count matrix with each origin or destination at the stop level. The optimization models to be presented identify each trip segment as either a singleton trip or a trip leg. The O-D matrix can be calculated elementwise and its $(i,j)^{th}$ entry represents the flow estimation from transit stop $i$ to $j$: Its value is the sum of the number of singleton trips that start in $i$ and end at $j$, and the number of multi-legged trips whose first leg starts in $i$ and last leg ends in $j$. It is also possible to construct the network O-D flows between origins and destinations at varying geographical resolutions based on the stop-level matrices. 

\section{Methodology}

This section presents optimization models to compute the aggregate O-D matrix. The optimization models do not compute the O-D matrix directly; instead they solve the Transfer Identification Problem (TIP) that identifies whether each trip segment is followed by a transfer, in which case the next trip segment is also identified. It is then simple to use the TIP solution to compute the aggregate O-D matrix. The optimization models choose the values of these decision variables in order to minimize the distance between the observed and estimated transfer probabilities at transit stops, subject to the transfer assumptions discussed earlier. 
This section presents three approaches for solving the TIP: A Quadratic Integer Program (QIP), a two-stage approach based on a continuous relaxation of the QIP and rounding, and an Integer Program (IP). 

\subsection{A QIP Model for the TIP}

Figure \ref{fig:opt} presents the QIP formulation for the TIP: the model specifies the setup and parameters, the decision variables, the objective function and the constraints, which are now discussed in detail.

\begin{figure*}[!t]
\begin{tabbing}
\tabrule\\
123\=123\=123\=12312312123341234\=12345678901234567\=123\=123\=123\=\kill
{\bf Data:} \\
\> $T$: set of recorded trips; \\
\> $C$: set of transit centers; \\
\> $p_{1,i}^*$: observed transfer probability at transit center $i$; \\
\> $p_{2}^*$: observed transfer probability at stops other than transit centers; \\
\> $t$: maximal transfer interval time in minutes; \\
\> $d$: maximal transfer walking distance in miles; \\
\> $\forall j \in T:$ \\
\>\> $l_j$: route of trip $j$; \\
\>\> $b_j$: boarding stop of trip $j$; \\
\>\> $a_j$: alighting stop of trip $j$; \\
\>\> $s_j$: boarding time of trip $j$; \\
\>\> $t_j$: alighting time of trip $j$; \\
\>\> $T_j := \{k \in T \mid  l_j \neq l_k, ~ dist(a_j, b_k) < d, ~ 0 < s_k - t_j < t \} $, set of possible transfers from $j$. \\

{\bf Variables:} \\
\> $x_j \in \{0,1\}$ \>\>\> $(j \in T)$ \> --- $x_j = 1$ if $j$ is a first leg; \\
\> $y_{j,k} \in \{0,1\}$ \>\>\> $(j \in T, k \in T_j)$ \> --- $y_{j,k} = 1$ if trip segment $j$ transfers to $k$; \\
\> $p_{1,i} \in [0,1]$ \>\>\> $(i \in C)$ \> --- calculated transfer probability at $i$; \\
\> $p_{2} \in [0,1]$ \>\>\> --- calculated transfer probability at stops other than transit centers. \\

{\bf Objective:} \\
\> minimize  ${\displaystyle \sum_{i \in C} \; (p_{1,i} - p_{1,i}^*)^2 + (p_{2} - p_{2}^*)^2}$ \\

{\bf Constraints:} \\

\> ${\displaystyle \sum_{k \in T_j} y_{j,k} = x_j}$ \>\>\> $(j \in T)$ \>\> (0.1) \\

\> ${\displaystyle \sum_{j \in T \mid k \in T_j}y_{j,k} \leq 1 - x_k}$ \>\>\> $(k \in T)$ \>\> (0.2) \\

\> ${\displaystyle p_{1,i} = \frac{\sum_{j \in T \mid a_j = i} \ x_j}{\mid \{j \in T \mid a_j = i \} \mid }}$ \>\>\> $(i \in C)$ \>\> (0.3) \\

\> ${\displaystyle p_{2} = \frac{\sum_{j \in T \mid a_j \notin C} \ x_j}{\mid \{j \in T \mid a_j \notin C \} \mid}}$ \>\>\>\>\> (0.4) \\

\tabrule
\end{tabbing}
\caption{The QIP Model for the TIP.}
\label{fig:opt}
\end{figure*}

The model is defined over the set $T$ of trip segments and the set $C$ of transit centers. Each trip segment $j \in T$ is characterized by its boarding stop $b_j$, alighting stop $a_j$, boarding time $s_j$, alighting time $t_j$, and bus line $l_j$. Two sets of decision variables are associated with each segment: Binary variable $x_j$ is 1 if and only if trip segment $j$ has a transfer and binary variable $y_{j,k}$ is 1 if and only if segment $j$ transfers to segment $k$. Recall that a transfer between two segments $j$ and $k$ has positive probability if and only if they are on different routes, the maximum walking distance and transfer time constraints are satisfied, and trip segment $k$ starts after $j$ ends. These constraints can be expressed as
$$
l_z \neq l_c \ \& \  dist(a_j, b_k) < d \ \& \  0 < s_k - t_j < t
$$

\noindent where the function $dist(\cdot,\cdot)$ is a metric (e.g., geodesic, Euclidean, or Manhattan distances), $d$ and $t$ are parameters chosen to denote the maximum allowed walking distance and transfer time. For each segment $j$, the set of feasible second legs from $j$ is denoted by $T_j$. 

Constraint (0.1) specifies that there is a transfer after trip segment $j$ (i.e. $x_j = 1$) if and only if there is exactly one transfer from this trip to one of its feasible successors (i.e., exactly one of $y_{j,k} = 1$ for $k \in T_j$). Constraint (0.2) states the one-transfer assumption: If segment $j$ transfers to segment $k$ (i.e., $y_{j,k} = 1$), then $k$ cannot make any subsequent transfers (i.e., $x_k = 0$). Conversely, if segment $k$ has a transfer (i.e., $x_k = 0$), then any feasible prior leg $j$ (i.e., all $j ~\textrm{such that}~ k \in T_j$) cannot transfer to $k$ (i.e., $y_{j,k} = 0$). Constraints (0.3) and (0.4) respectively compute the estimated transfer probabilities at each transit center and all other stops, where $\sum_{j \in T \mid a_j = i} \ x_j$ represents the number of first legs transferring at transit center $i$ and 
$\mid \{j \in T \mid a_j = i \} \mid$ counts the total number of segments ending at transit center $i$. Recall that the goal of the TIP is to select potential transfers for each trip segment so that the aggregate-level transfer probabilities given by Constraints (0.3) and (0.4) are as close as possible to the observed transfer probabilities. In the QIP, closeness is measured with an L2-norm. The sets of first legs, second legs, and singleton trips are denoted respectively by $T_1$, $T_2$ and $T_s$, which form a partition of the set $T$. They are defined explicitly as follows, 
\begin{align*}
    T_1 &:= \{j \in T \mid x_j = 1\} \\
    T_2 &:=\{k \in T \mid y_{j,k} = 1 ~\textrm{for some}~ j \in T\} \\
    T_s &:= T \setminus (T_1 \cup T_2)
\end{align*}

The O-D matrix can then be estimated elementwise from the QIP solution. Each entry $(i, i')$ of the O-D matrix records the expected number of trips from transit stop $i$ to $i'$ and is a sum of two components: The number of singleton trips starting at stop $i$ and ending at $i'$ and the number of two-legged trips whose first leg starts at stop $i$ and second leg ends at $i'$. Let $\textbf{1}_{ \{ \cdot \} }$ be an indicator function which equals to 1 when the statement is true and 0 otherwise. Then the O-D matrix can be computed as follows:
\begin{align*}
    OD_{i,i'} = &\sum_{j \in T_s} \textbf{1}_{\{b_j=i,~a_j=i'\}} \\ + &\sum_{j \in T_1} \sum_{k \in (T_j \cap T_2)} y_{j,k} \textbf{1}_{\{b_j=i,~a_k=i'\}}.
\end{align*}

The proposed QIP formulation is applied to the case study to be introduced in Section IV and shown to have severe scalability issues: It cannot be solved by Gurobi \cite{gurobi}, a state-of-the-art commercial optimization solver, within 24 hours. In general, the QIP formulation does not have guaranteed tractability due to the (potentially) quadratic number of variables and the large size of the data sets, even for small-size cities.

\subsection{Rounding the Continuous QIP Relaxation}

This section explores a scalable two-stage approach which consists of (1) solving the continuous relaxation of the QIP and (2) rounding the solutions of continuous relaxation to derive a feasible binary substitute. The continuous relaxation relaxes the domain of the variables from the set $\{0,1\}$ to the interval $[0,1]$, producing a Convex Quadratic Program (QCP), which can be solved efficiently. The QCP solution for variables $x_j$ and $y_{j, k}$ $(j \in T, k \in T_j)$ now assigns values in the range $[0,1]$ to the decision variables which can thus be interpreted as the probability of having a transfer after segment $j$ and the probability that segment $k$ be the second leg of that transfer. The QCP relaxation is depicted in Figure \ref{fig:opt1} and mimics the QIP.

\begin{figure*}[!ht]
\begin{tabbing}
\tabrule\\
123\=123\=123\=12312312123341234\=12345678901234567\=123\=123\=123\=\kill

{\bf Variables:} \\
\> $x_j \in [0,1]$ \>\>\> $(j \in T)$ \> --- probability of $j$ being a first leg; \\
\> $y_{j,k} \in [0,1]$ \>\>\> $(j \in T, k \in T_j)$ \> --- probability of $j$ transferring to $k$; \\
\> $p_{1,i} \in [0,1]$ \>\>\> $(i \in C)$ \> --- calculated transfer probability at $i$; \\
\> $p_{2} \in [0,1]$ \>\>\> --- calculated transfer probability at stops other than transit centers. \\

{\bf Objective:} \\
\> minimize  ${\displaystyle \sum_{i \in C} \; (p_{1,i} - p_{1,i}^*)^2 + (p_{2} - p_{2}^*)^2}$ \\

{\bf Constraints:} \\

\> ${\displaystyle \sum_{k \in T_j} y_{j,k} = x_j}$ \>\>\> $(j \in T)$ \>\> (1.1) \\

\> ${\displaystyle \sum_{j \in T \mid k \in T_j}y_{j,k} \leq 1 - x_k}$ \>\>\> $(k \in T)$ \>\> (1.2) \\

\> ${\displaystyle p_{1,i} = \frac{\sum_{j \in T \mid a_j = i} \ x_j}{\mid \{j \in T \mid a_j = i \} \mid }}$ \>\>\> $(i \in C)$ \>\> (1.3) \\

\> ${\displaystyle p_{2} = \frac{\sum_{j \in T \mid a_j \notin C} \ x_j}{\mid \{j \in T \mid a_j \notin C \} \mid}}$ \>\>\>\>\> (1.4) \\

\tabrule
\end{tabbing}
\caption{The QCP Relaxation of the TIP.}
\label{fig:opt1}
\end{figure*}

To obtain an aggregate O-D matrix, it is necessary to round the variables and assign them binary values. The second stage is based on a feasible rounding strategy that proceeds as follows. First, the segments likely to have a transfer are rounded to 1 by choosing a threshold $x^* \in [0,1]$ and selecting those variables whose value in the QCP relaxation exceeds the threshold, i.e.,
$$
\hat{x}_j = \begin{cases} 1, &if~x_j \geq x^* \\ 0, &otherwise. \end{cases}
$$

\noindent
The threshold can be obtained by rounding up the variables $x_j~(j \in T)$ with $n$ largest probability of having a transfer to 1 and rounding down the rest to 0. For notation simplicity, define the observed transfers: $\delta_{1,i}$ for transit center $i$,  $\delta_2$ for other non-transit stops and $n$ for the total number of transfers based on the observed transfer probabilities as follows,
$$
\begin{array}{ll}
\delta_{1,i} &= p^*_{1,i} \times |\{j \in T \mid a_j = i \}| \\
\delta_2 &=  p^*_2 \times |\{j \in T \mid a_j \notin C \}| \\
n &= \sum_{i \in C} \delta_{1,i} + \delta_2.
\end{array}
$$
The threshold $x^*$ can be derived as below,
$$ x^* = \min_{j \in T} x_j ~\textrm{s.t.}~ |\{j' \in T|x_{j'} \geq x_j \}| \leq n ,$$
and the set of first legs consists of all trip segments indexed by $j \in T$ with the corresponding variable $x_j \geq x^*$, that is,
$$ T_1 := \{ j \in T | \hat{x}_j = 1 \} = \{ j \in T | x_j \geq x^* \}. $$



It remains to determine the set of second legs and the set of singleton trips. The likelihood of a segment to be  a second leg can be approximated by summing the transfer probabilities from all of its possible first legs, i.e., $\mathcal{L}_k := \sum_{j \in T_1} y_{j,k} ~(k \in T \setminus T_1)$. Note that, since the transfer probabilities from different first legs are computed in different probability spaces, their sums do not directly translate to a probabilistic interpretation. However, it is still a sensible measure for identifying the set of segments most likely to be second legs. Recall that the one-transfer assumption requires that the number of second legs equals the number of first legs and that $n$ as defined earlier denotes the total observed transfers: The set of second legs can now be defined by those segments $k$ whose measure $\mathcal{L}_k$ is among the $n$ largest values , i.e.,
$$
T_2 := \{k \in T \setminus T_1 \mid ~~ \big \lvert \{k' \in T\setminus T_1 \mid \mathcal{L}_{k'} \geq \mathcal{L}_k\} \big \rvert \leq n\}.
$$


Once the trips are grouped as first-leg, second-leg, or singleton-trip, the probabilities of transferring from a first-leg trip $j$ to any feasible second leg $k \in (T_j \cap T_2)$ are normalized to
produce, for each first leg, a well-defined probability distribution over its feasible second legs, i.e., $\forall j \in T_1$, $\sum_{k \in T_j \cap T_2} p_{j,k} = 1$ and $p_{j,k} \geq 0 ~~(k \in T_j \cap T_2)$. The normalized probabilities are calculated as,
$$ p_{j,k} = \frac{y_{j,k}}{\sum_{k' \in T_j \cap T_2} y_{j,k'}}~~ \forall k \in T_j \cap T_2, \forall j \in T_1. $$
The O-D matrix can then be estimated elementwise as a probability-weighted sum of trips. The expected number of trips from transit stop $i$ to $i'$ is the sum of two components: The count of singleton trips starting at stop $i$ and ending at $i'$, and the sum of (normalized) probabilities for all feasible two-legged trips whose first leg starts at stop $i$ and the second leg ends at $i'$, i.e.,
\begin{align*}
    OD_{i,i'} = &\sum_{j \in T_s} \textbf{1}_{\{b_j=i,~a_j=i'\}} \\ + &\sum_{j \in T_1} \sum_{k \in (T_j \cap T_2) } p_{j,k} \textbf{1}_{\{b_j=i,~a_k=i'\}}.
\end{align*}

\subsection{The Integer Programming (IP) Model}

This section proposes a third approach based on Integer Programming (IP). The key idea behind the IP model is to replace the L2-norm by a L1-norm and to reason about the \textit{observed transfers} instead of the \textit{observed transfer probabilities}. Recall that the observed transfers are defined based on the observed transfer probabilities as $\delta_{1,i}$ for each transit center $i \in C$ and $\delta_2$ for other non-transit stops.
Figure \ref{fig:opt2} describes the resulting IP model. The objective function minimizes the absolute differences between the observed and estimated numbers of transfers. The logical constraints are the same as in the QIP, but there is no need to reason about transfer probabilities. The aggregate O-D matrix can then be computed from the optimal solution as for the QIP.

\begin{figure*}[!ht]
\begin{tabbing}
\tabrule\\
123\=123\=123\=12312312123341234\=12345678901234567\=123\=123\=123\=\kill
{\bf Data:} \\
\> $\delta_{1,i} \in \mathbb{N} $ \>\>\> $(i \in C)$ \> --- observed transfers at transit center $i$; \\
\> $\delta_{2} \in \mathbb{N}$ \>\>\>\> --- observed transfers at other stops. \\ 
{\bf Variables:} \\
\> $x_j \in \{0,1\}$ \>\>\> $(j \in T)$ \> --- $x_j = 1$ if trip $j$ transfers; \\
\> $y_{j,k} \in \{0,1\}$ \>\>\> $(j \in T, k \in T_j)$ \> --- $y_{j,k} = 1$ if trip $j$ transfers to $k$; \\

{\bf Objective:} \\
\> minimize  ${\displaystyle \sum_{i \in C} \| \delta_{1,i} - \sum_{j \in T \mid a_j = i} x_j \|_1 +  \| \delta_2 - \sum_{j \in T \mid a_j \notin C} x_j \|_1}$ \\

{\bf Constraints:} \\

\> ${\displaystyle \sum_{k \in T_j} y_{j,k} = x_j}$ \>\>\>\> $(j \in T)$ \>\>\> (2.1) \\

\> ${\displaystyle \sum_{j \in T \mid k \in T_j}y_{j,k} \leq 1 - x_k}$ \>\>\>\> $(k \in T)$ \>\>\> (2.2) \\
\tabrule
\end{tabbing}
\caption{The Integer Programming Model for the TIP.}
\label{fig:opt2}
\end{figure*}

\section{Case Study}
This section applies the proposed methodology to a case study for the broader Ann Arbor--Ypsilanti region in Michigan and validates the O-D flow matrices estimated from both the rounded QCP and the IP.

\begin{figure*}[!t]
    \includegraphics[width=16cm]{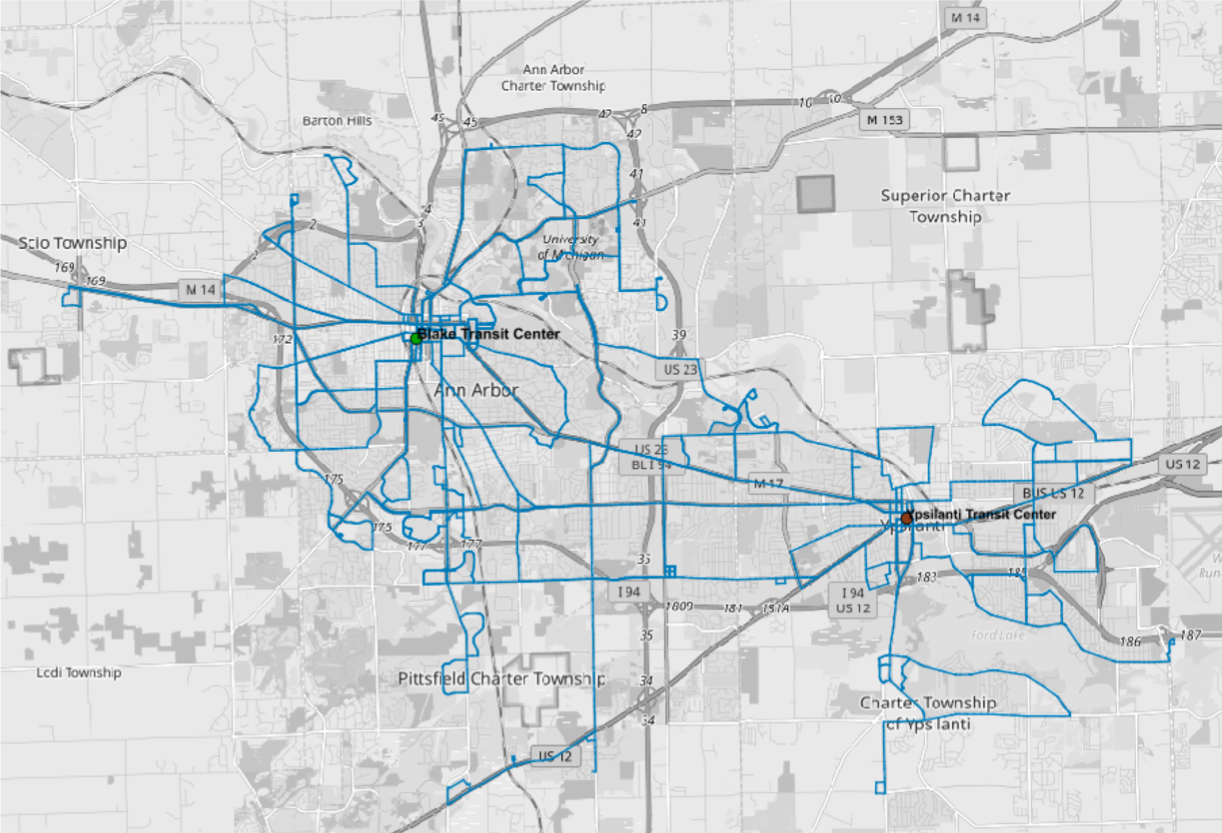}
    \centering
    \caption{Transit Network Operated by AAATA}
    \label{fig:aaata}
\end{figure*}

\subsection{Data Description} The data in this study was provided by the Ann Arbor Area Transportation Authority (AAATA), consisting of boarding-only smart-card transactions. The data was collected from the transit network operated by AAATA as depicted in Figure \ref{fig:aaata}, where a total of 58 inbound or outbound routes (denoted by blue lines in Figure \ref{fig:aaata}) connect 1,232 stops including two transit centers, namely Blake Transit Center and Ypsilanti Transit Center. The bus schedule and stop location information were extracted from the GTFS data \cite{TheRide2019}. 
We conducted two experiments using the Go!pass data and the Period Pass data, respectively. Go!pass is purchased by local businesses at downtown Ann Arbor as a benefit for their employees who can get unlimited bus services provided by AAATA. Period Pass, on the other hand, allows the holders to take unlimited rides on any bus routes within the specified period and offers discounted prices for senior and students \cite{TheRide2019}. There are overall 32,840 transactions for Go!Pass and 43,660 for Period Pass from all weekdays in October, 2017.

\subsection{Experimental Settings}  

In this study, we use the smart-card data to validate the proposed methodology. To be specific, we apply the trip-chaining method on the boarding-only smart-card data to infer alighting stops, identify transfers, and produce a transit O-D matrix to evaluate and validate our optimization models \cite{trepanier2007individual,munizaga2012estimation}. Note that we assume the ground truth is given by the trip-chaining method on the smart-card data.
The trip-chaining method exploits unique IDs for passengers to link their consecutive transit trip segments. It relies on two major assumptions: 1) the alighting point of each trip is within walking distance of the consequent boarding location (usually assumed to be in the range of 400 m to 2,000 m); 2) the alighting point of the last boarding of the day for a passenger is adjacent to his/her first boarding stop of the same day. In addition, researchers typically assume a time threshold (e.g., 30 minutes) to identify transfer activities: A passenger is assumed to take a transfer if the interval between the alighting time and the subsequent boarding is less than the specified threshold \cite{devillaine2012detection}. For the trip-chaining benchmark, we assume the maximum walking distance is a quarter mile (402 meters), transfer time threshold is 30 min, and the last destination is assumed to be the closest stop to the first origin as suggested in \cite{alsger2016validating}. 

Recall that the optimization models require all route-level time-stamped O-D matrices collected from a transit network as its input, which can be readily calculated from APC data using the IPF techniques. Therefore, to ensure the validity of the comparison, we process the same smart-card data to generate the route-level O-D matrices by directly aggregating the inferred route-level boarding-alighting pairs for each transaction into time-stamped route-level O-D flows, which will serve as the input for our optimization models. To be consistent with the settings of the benchmark, the optimization models also assume the same maximum walking distance (402 meters) and transfer time thresholds (30 min). 





There are two transit centers in the region, i.e., one in downtown Ann Arbor, the other in downtown Ypsilanti. We assume that the Blake Transit Center (in Ann Arbor) takes $i = 1$ and the Ypsilanti Transit Center takes $i = 2$.
The ground-truth transfer rates computed from the smart-card analysis are: For Go!Pass, 
$$P_{11}^* = 0.232, ~P_{12}^* = 0.591 ~\textrm{and}~ P_2^* = 0.062;$$
for Period Pass, 
$$P_{11}^* = 0.588, ~P_{12}^* = 0.554 ~\textrm{and}~ P_2^* = 0.143.$$ 

\subsection{Geographical Resolutions for Model Evaluation}

The estimated O-D matrices by our models are evaluated at various geographical resolutions, i.e., the stop level, the TAZ level (which is widely adopted for transportation planning), and the Transit Analysis Clusters (TACs). TACs are self-defined zones, obtained by using Hierarchical Clustering Analysis (HCA) \cite{HCA2013}. In HCA, each stop is initially assigned to a cluster on its own, and at each step, two most similar clusters (based on distance) are joined until a single cluster is left. The end result is a giant cluster organized in a tree structure. To obtain TACs at different geographical resolutions, it suffices to cut the tree at a given height: The resultant clusters have pairwise distances approximately at that height which is chosen to capture the desired distance threshold.

There are three major reasons behind this choice of TACs through HCA. First, riders may choose different origin and destination stops via different bus routes in a day, so aggregating the travel demand of spatially close stops is critical when constructing O-D matrices \cite{luo2017constructing}. Second, the stop-to-stop O-D matrix does not directly reveal the travel demand pattern like a zone-to-zone O-D matrix which shows the equilibrium between demand and supply in the transit system and measures the access or egress times for transit trips \cite{tamblay2016zonal}. Third, traditional TAZs have boundaries on the streets, where most bus stops are located. Hence the TAZs boundaries might create artificial divisions between closely related stops and influence the estimation. We will present the evaluation metrics at stop level, TAZ level and TAC level with varying radius to provide a more comprehensive examination of model performance.  

\subsection{Experimental Results}

The R-squared metric is used to evaluate the accuracy of the O-D flow estimates as in previous studies \cite{Tavassoli2016HowCT} and \cite{EconometricsInTransportation}. Given a ground-truth square matrix $OD^{*}$ and its estimation $OD$, both of which with size of $n \times n$, the R-squared metric is defined as
$$
\theta =  1 - \frac{\sum_{i,j=1}^n (OD^{*}_{ij} - OD_{ij})^2 }{\sum_{i,j=1}^n (OD^{*}_{ij} - \overline{OD^{*}})^2 }.
$$

\noindent where $\overline{OD^{*}}$ denotes the average of all entries of $OD^{*}$. The R-squared metric can be interpreted as the percentage of total variability (in terms of sum of squared errors) that can be explained by the estimated matrix, compared to a mean model. A higher R-squared value generally indicates a better estimation. 

As discussed, the models are evaluated at different geographical resolutions: The stop level, the TAZ level, and the TAC level with radius ranging from 0.5 miles to 2 miles in increments of 0.5 miles. The stop level can also be seen as a TAC level with radius of 0 mile. The results at the TAZ level are now summarized. For the Go!Pass data, the R-squared evaluated against the ground-truths inferred from trip-chaining methods is equal to 88.71\% for the rounded QCP model and 95.57\% for the IP model. For the Period Pass data, the R-squared reaches 67.15\% for rounded QCP model and 85.06\% for the IP model. Figure 5 illustrates the R-Squared evaluated using both Go!Pass and Period Pass data for the radius of TAC ranging from 0 mile to 2 miles.




\begin{figure*}[t!]
    \includegraphics[width=16cm]{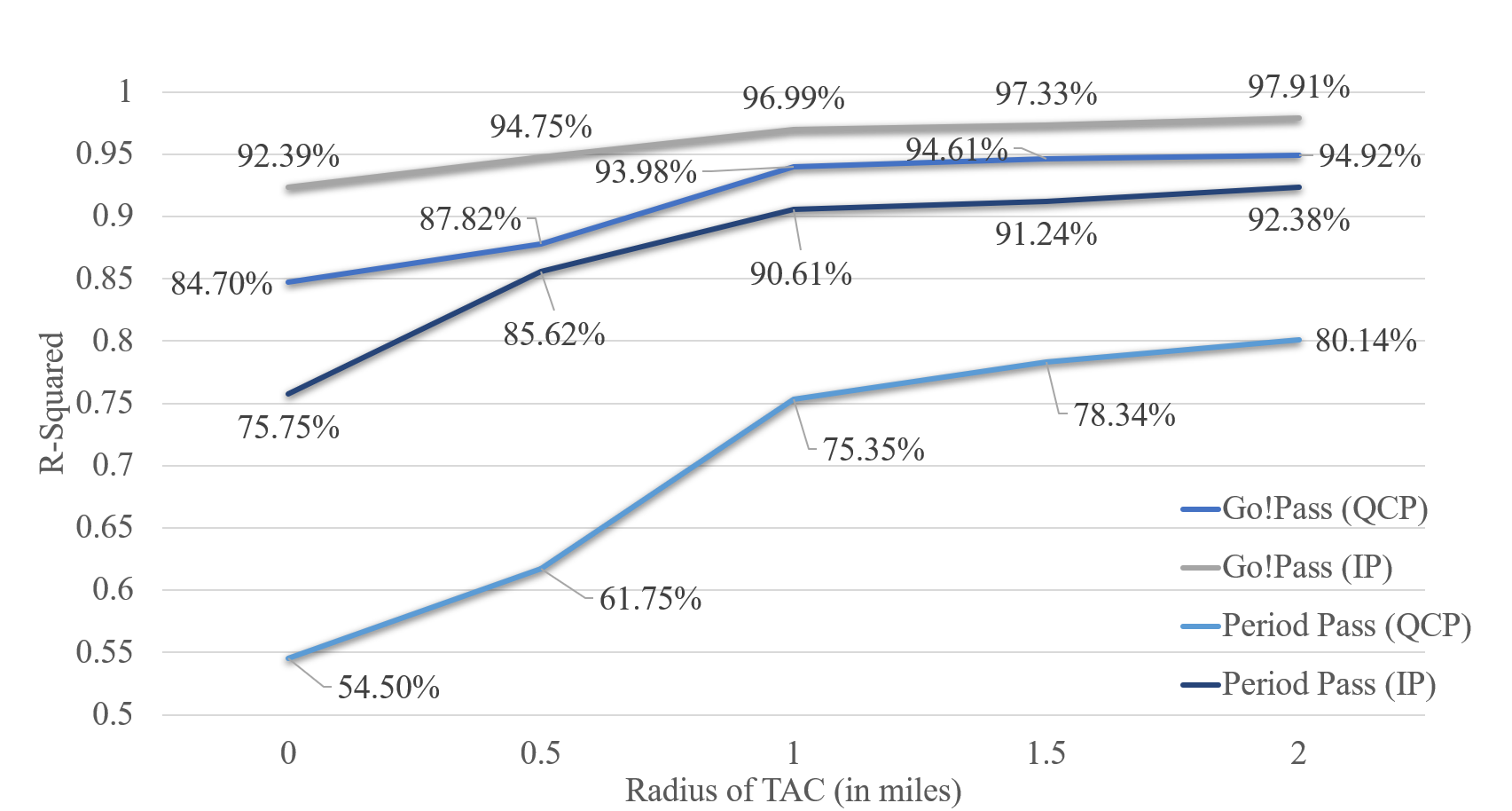}
    \centering
    \caption{R-Squared Values at Varying Radius of TACs.}
    \label{fig:rsq}
\end{figure*}

Both approaches can be solved efficiently on a personal computer of i7-7500U CPU with Gurobi 8.0 through Python 3.6. The rounded QCP approach takes 121.78 seconds on Go!Pass data (with 32,840 transactions) and 156.96 seconds on Period Pass data (with 43,660 transactions); the IP approach solves faster with a computational time of 2.34 seconds on Go!Pass data and 5.97 seconds on Period Pass data.

There are three key general observations. First, the IP approach performs significantly better than rounded QCP approach on both data sets. Second, the predictive power of both approaches improves with increasing TAC radius but the rate of improvement typically decreases as the clusters expand in their radius. The Rounded QCP approach on Go!Pass increases monotonically from 84.70\% at 0-mile radius to 94.92\% at 2-mile radius. A similar trend can be observed from the IP approach on the Go!Pass data where $\theta$ increases monotonically from 92.39\% at 0-mile radius to 97.91\% at 2-mile radius. For both models, the same behavior is also observed on the Period Pass data. Furthermore, note that the accuracy of 1-mile TAC is an important turning point as shown in Figure \ref{fig:rsq}---the improvement rate below the 1-mile threshold is significantly larger than above it.
Third, the model performance on the Period Pass data is worse than the Go!Pass data, which might result from the irregular travel behavior of Period Pass users. Recall that the Go!Pass is a special pass purchased by companies located in downtown areas for their daily commuting employers, who tend to have more consistent and predictable commuting patterns during the weekdays. Therefore, for the Go!Pass data, the daily variation in transfers at each stop may be less and can be better described by a single observed long-term transfer rate. In comparison, Period Pass targets a more general audience and offers special discounts for the seniors and students who tend to take more spontaneous and thus less predictable trip chains. It is also observed that the true transfer rates for some days in the Period Pass data deviate dramatically from the long-term observation used in the optimization models. Also, the transfer rates for Period Pass is much higher than Go!Pass, which might lead to more feasible transfers for each first-leg trips and add to the modeling complexity. These discussions might indicate an inherent limitation in modeling transfer identification with inadequate information. Without unique identifier of passengers, it is expected that when dealing with extremely high transfer probabilities (such as 50\% or above), the aggregate-level variability is difficult to capture in detail.

\section{Conclusion}

This paper presented optimization models to estimate the transit network O-D flow matrix based on time-stamped and location-based boarding and alighting counts, and observed or estimated proportions of transferring passengers at each transit center and other non-transit centers. It proposed a QIP approach, a two-stage approach based on a QCP relaxation of the QIP and a feasible rounding procedure, and an IP model that replaces the L2-norm of the QIP by a L1-norm. While the QIP is not tractable for real data sets, the QCP and IP approaches can be solved efficiently for the transit data provided by AAATA. Moreover, the IP model is superior to the QCP approach in terms of accuracy.

When measuring against the ground-truth as calculated from trip-chaining methods using the R-squared metric, the IP model can achieve up to 95.59\% at the TAZ level and 96.99\% at the 1-mile self-defined TAC level for the Go!Pass data (which exhibits more consistent travel behavior) and 85.06\% at the TAZ level and 90.61\% at the 1-mile TAC level for the Period Pass data (which exhibits more irregular travel patterns). There is also a clear improvement in predictive accuracy with lower spatial resolution. The results suggest that the IP model can produce accurate estimation for applications requiring varying levels of spatial resolutions. In particular, the IP model can meet the needs for tasks ranging from predicting O-D flow among bus stops to constructing a zone-level transit-trip O-D matrix to inform future transit planning. The results indicate the IP model is especially promising for cases with moderate or relatively low transfer rates and for populations with consistent transfer patterns. It is because the observed transfer probabilities are used as benchmarks, against which the deviation of estimated transfer probabilities is minimized. Therefore, the capability of observed transfer probabilities to capture a consistent pattern for transfer activities is a critical factor for accurate modeling. It is recommended to apply the transfer probabilities observed or surveyed during a relatively long period (i.e., monthly or yearly), which reduces the potential inaccuracy brought by daily variations. If such information is absent, the expert judgment from transit operators could be used instead. 

This study can be further developed from the following perspectives. First, the IP model directly applies the parameters for behavioral assumptions as suggested in \cite{alsger2016validating} for their case study on South-East Queensland public transport network in Australia. However, the parameters for such behavioral assumptions might differ across case studies due to differences in transit systems, built environments and socio-demographics of the regions under analysis. More comprehensive studies to validate the proposed methodology would welcome experimental results on more case studies or a sensitivity analysis on the choice of maximally allowed distance and transfer time. 
Second, the current IP formulation can be easily extended to account for multiple transfers. Future work can verify the effectiveness of integer programming modeling multiple transfers. Also, the optimization models for transfer identification may have many symmetric solutions, leading to large deviations when transfer probabilities are high. As a result, it is important to clearly identify these equivalent solutions and conduct closer case-dependent analysis to obtain a more accurate prediction of the O-D pairs. This is a key direction for future research.

The methodology proposed in this study mainly serves to extend the current analysis of AVL/APC data and produce a network-level O-D matrix to inform transportation planning. Our models will also be suitable for analyzing smart-card data or Automatic Farebox Collection (AFC) data with hidden unique identifiable information. Recently, as stated in \cite{pelletier2011smart}, the use of smart-card data has raised privacy concerns. One major problem is the vulnerability of the central database which stores smart-card transactions and user information, especially when the data is used for multiple purposes and accessible by multiple groups. Withholding unique ID information when releasing the data to third parties could significantly reduce the risk of private information disclosure. However, lacking unique ID information would prevent the use of the trip-chaining methods. Therefore, the transfer identification model provides a tangible tool for estimating travel demand for such data at an aggregate level.


\section*{Acknowledgment}
This research is funded by the Michigan Institute of Data Science (MIDAS) and by Grant 7F-30154 from the Department of Energy. The authors would like to thank Forest Yang from the AAATA for his assistance in providing the data. Findings presented in this paper do not necessarily represent the views of the funding agencies.

\ifCLASSOPTIONcaptionsoff
  \newpage
\fi



\bibliographystyle{IEEEtran}
\bibliography{main.bib}
\end{document}